\documentclass[12pt]{amsart}
\usepackage{amssymb,amsmath,mathrsfs}
\usepackage[ps2pdf,colorlinks=true,urlcolor=blue,
citecolor=red,linkcolor=blue,linktocpage,pdfpagelabels,bookmarksnumbered,bookmarksopen]{hyperref}
\usepackage{epsfig,graphicx,color,mathrsfs}
\usepackage[english]{babel}

\newcommand{\N}{{\mathbb N}}

\newcommand{\D}{{\mathbb D}}
\newcommand{\R}{{\mathbb R}}

\newcommand{\be}{\begin{equation}}
\newcommand{\ee}{\end{equation}}

\newcommand{\eps}{\varepsilon}

\newcommand{\calm}{{\mathcal M}}

\numberwithin{equation}{section}
\newtheorem{theorem}{Theorem}[section]
\newtheorem{proposition}[theorem]{Proposition}

\newtheorem{lemma}[theorem]{Lemma}

\theoremstyle{definition}
\newtheorem{remark}[theorem]{Remark}

\newcommand{\brm}{\begin{remark}\rm}
\newcommand{\erm}{\end{remark}}
\newcommand{\brms}{\begin{remark}\rm}
\newcommand{\erms}{\end{remark}}
\newcommand{\bte}{\begin{theorem}}
\newcommand{\ete}{\end{theorem}}
\newcommand{\bpr}{\begin{proposition}}
\newcommand{\epr}{\end{proposition}}
\newcommand{\ble}{\begin{lemma}}
\newcommand{\ele}{\end{lemma}}
\newcommand{\beq}{\begin{equation}}
\newcommand{\eeq}{\end{equation}}
\newcommand{\bdm}{\begin{displaymath}}
\newcommand{\edm}{\end{displaymath}}
\numberwithin{equation}{section}

\newcommand{\bos}{\begin{remark}\rm}
\newcommand{\eos}{\end{remark}}

\newcommand{\ben}{\begin{enumerate}}
\newcommand{\een}{\end{enumerate}}

\title[Minimization under constraints: the added
mass technique]{An approach to minimization under \\ constraint: the
added mass technique}

\author[Louis Jeanjean]{Louis Jeanjean$^*$}
\address{Laboratoire de Math\'ematiques (UMR 6623)
\newline\indent
Universit\'{e} de Franche-Comt\'{e}
\newline\indent
16, Route de Gray 25030 Besan\c{c}on Cedex, France}
\email{louis.jeanjean@univ-fcomte.fr}

\thanks{$^*$Laboratoire de Math\'ematiques (UMR 6623),
Universit\'{e} de Franche-Comt\'{e},
16, Route de Gray 25030, Besan\c{c}on Cedex, France.
E-mail: {\em louis.jeanjean@univ-fcomte.fr}}

\author[Marco Squassina]{Marco Squassina$^\dagger$}
\address{Dipartimento di Informatica
\newline\indent
Universit\`a degli Studi di Verona
\newline\indent
C\'a Vignal 2, Strada Le Grazie 15, 37134 Verona, Italy}
\email{marco.squassina@univr.it}

\thanks{$^\dagger$Dipartimento di Informatica,
Universit\`a di Verona,
C\`a Vignal 2, Strada Le Grazie 15, I-37134 Verona, Italy.
E-mail: {\em marco.squassina@univr.it}}

\thanks{The second author was partially supported by the
Italian PRIN Research Project 2007 {\em Metodi variazionali e topologici
nello studio di fenomeni non lineari}}

\begin{document}
\subjclass[2000]{35J40; 58E05}

\keywords{Constrained minimization problems,
concentration compactness, quasi-linear elliptic
equations and systems}

\begin{abstract}
We present an approach to minimization under constraint. We
explore the connections of this technique with the general method
of Compactness by Concentration of P.L.\ Lions~\cite{lions1} and
present applications to some constrained semi-linear and quasi-linear
elliptic problems.
\end{abstract}
\maketitle



\section{Introduction}
In this paper we discuss an approach for the minimization of
functionals under a constraint and give some applications of it. We
start with a simple statement in order to illustrate our technique.
Let $H$ be a reflexive Banach function space on $\R^N$ ($N\geq 1$)
with value in $\R^m$ ($m\geq 1$) and let $J$, $G$ be functionals
defined on $H$ of the type
$$
J(u) = \int_{\R^N}j(x,u,|\nabla u|) dx,  \quad G(u) =
\int_{\R^N}g(u) dx,
$$
where $j(x,s,t)$ and $g(s)$ are real-valued functions defined on
$\R^N \times \R^m \times \R$ and $\R^m$ respectively. For a fixed
$c \in \R$, we consider the problem
\begin{eqnarray}
    \label{00.1}
\text{minimize $J$ on the functions $u\in H$ with $G(u) = c$}.
\end{eqnarray}
Setting
$$
m(c) =  \inf  \{ J(u): \text{$u\in H$ with $G(u) = c$}\},
$$
we have the following

\begin{proposition}
    \label{abstract-result}
Assume that $m(c) >- \infty$ and that there exists a minimizing
sequence $(u_n) \subset H$ such that
\begin{enumerate}
\item[(H0)] $(u_n) \subset H$ is bounded in $H$.
\item[(H1)] If $u_n \rightharpoonup u$ then
$$
J(u) \leq \liminf_{n \to \infty} J(u_n) \quad \mbox{ and } \quad  G(u) \leq c.
$$
\end{enumerate}
Then $m(c)$ is reached if in addition
\begin{enumerate}
\item[(H2)] There exists $v \in H$ such that
$$
G(u+v) = c \quad \mbox{ and } \quad J(u+v) \leq J(u).
$$
\end{enumerate}
\end{proposition}

\begin{proof}
Let $(u_n) \subset H$ satisfy $(H0)$. Then $(u_n) \subset H$ is
bounded and we can assume that, up to a subsequence, $u_n
\rightharpoonup u$ in $H$, for some $u \in H$. Then by (H1) we get
that
$$
J(u) \leq \liminf_{n \to \infty}J(u_n) = m(c)
$$
with $G(u) \leq c$. If $G(u) = c$ we are done (and condition (H2) holds with
$v=0$). If $G(u)<c$  by (H2) there exists a $v \in H$ such that
$G(u+v) = c$ and $J(u+v) \leq J(u) \leq m(c)$. If  $J(u+v) < m(c)$
this contradicts the definition of $m(c)$. Hence $J(u+v) = m(c)$,
so that $(u+v)$ is a minimizer for $m(c)$.
\end{proof}

Of course, assumption (H0) is necessary to study the minimization
problem~\eqref{00.1}. The fact that assumption (H1) holds, for at
least a bounded minimizing sequence, is more restrictive and
somehow defines the class of minimization problems under study.
The third assumption (H2) is clearly necessary for $m(c)$ to be
reached. Indeed if $u_0$ is a minimizer of $m(c)$ then taking
$v=0$ we have $G(u_0+v)= c$ and $J(u_0 +v) = J(u_0) = m(c)$. We use
assumption (H2) in the following way. Assuming, by contradiction,
that the weak limit $u \in H$ obtained in (H1) is not a minimizer
we construct a $v \in H$ such that $G(u+v) = c$ and $J(u+v) <J(u)
\leq m(c)$. Namely, checking (H2) relies on the possibility to
``add mass", that is to increase $c$, while strictly decreasing
the value of the functional $J$.
\medskip

In order to motivate the introduction of
Proposition~\ref{abstract-result} we first state the following
result. It is a special case of Proposition~\ref{abstract-result},
which is also useful by itself.

\begin{proposition}\label{Proposition2}
Assume that conditions (H0)-(H1) hold and that the function
$\lambda \in \R \mapsto m(\lambda)$ is strictly decreasing. Then,
for any fixed $c \in\R^+$, the value $m(c)$ is reached.
\end{proposition}

\begin{proof}
Let $c\in\R$ be fixed. By (H0) there exists a bounded minimizing
sequence $(u_n) \subset H$ and we can assume that $u_n
\rightharpoonup u$ in $H$ as $n\to\infty$. From (H1) we get that
$J(u) \leq m(c)$. Thus necessarily we obtain $m(G(u)) \leq m(c)$
and so, if it was $G(u) <c$, we would get a contradiction with the
assumption that the map $\lambda \mapsto m(\lambda)$ is strictly
decreasing.
\end{proof}

Over the last twenty five years the Compactness by Concentration
of P.L.\ Lions~\cite{lions1} has had a deep influence on the
problem of minimizing a functional under a given constraint. Let us
assume, for the moment, that we can define a problem at infinity
associated to~\eqref{00.1}. The limit of $j(x,u,|\nabla u|)$ as $|x|\to\infty$ is
denoted $j_{\infty}(u, |\nabla u|)$ and, accordingly, we define
$$
J_{\infty}(u) = \int_{\R^N} j_{\infty}(u, |\nabla u|)dx
$$
and
$$
m_{\infty}(c) =  \inf \{ J_{\infty}(u): \text{$u\in H$ with $G(u) = c$}\}.
$$
In~\cite{lions1} it is shown that all minimizing sequences
for~\eqref{00.1} are compact if, and only if, the following strict
inequality holds
\begin{equation} \label{strict-inequalities}
m(c) < m(\lambda) + m_{\infty}(c-\lambda), \quad \forall \lambda
\in [0,c[.
\end{equation}
The information that {\em all} minimizing sequences are compact is
essential in many situations, in particular when one deals with
orbital stability issues (see, for example,~\cite{cazelions}).
However if the issue is merely the existence of a minimizer one has
the freedom to choose {\em a particular} minimizing sequence. In
Propositions~\ref{abstract-result} and~\ref{Proposition2} we
exploit this fact and this allows us to treat cases which may not
satisfy condition~\eqref{strict-inequalities}.  In~\cite{lions1} it is
also heuristically explained (see pages 113-114) that the
corresponding large inequalities
\begin{equation}
    \label{large-inequalities}
m(c) \leq m(\lambda) + m_{\infty}(c-\lambda), \quad \forall c>0,
\quad \forall \lambda \in [0,c[
\end{equation}
are expected to hold under very weak assumptions. A direct
consequence of~\eqref{large-inequalities} is that, if
$m_{\infty}(d) <0$ for any $d \in [0,c[$, then the function
$\lambda \mapsto m(\lambda)$ is strictly decreasing. Thus we see, from
Proposition~\ref{Proposition2}, that in this case $m(c)$ is
reached just under (H0) and (H1). However in many situations the
condition $m_{\infty}(d) <0$ for any $d \in [0,c[$ is either difficult to
check or does not hold. On the contrary, proving that $m_{\infty}(d)
\leq 0$ for any $d \in [0,c[$, is often much easier. Note that,
following the heuristic discussion of~\cite{lions1}, we can then
still deduce that $\lambda \mapsto m(\lambda)$ is non increasing.
Knowing that the function $\lambda \mapsto m(\lambda)$ is non
increasing is often very useful to check assumption (H2) on
specific examples. Indeed, by applying the approach of
Proposition~\ref{abstract-result}, we can assume that there exists
a minimizing sequence $(u_n) \subset H$, $u_n \rightharpoonup u$ as $n\to\infty$,
for which
$$
J(u) \leq m(c),\quad \, \mbox{ with \,\,\,$G(u) \leq c$}.
$$
Then, if we can find a function $v \in H$ with $G(u) \leq G(v)
\leq c$ and $J(v) < J(u)$, we get a contradiction that proves that
$m(c)$ is reached.
There are also minimization problems which do not admit a
``problem at infinity" and thus where the approach
of~\cite{lions1} does not work. Also, in some cases, applying the
approach~\cite{lions1} leads to long proofs which could be shortened.
Ultimately, we point out that some of the ideas of this paper recently turned out
to be useful in the study of orbital stability for a class of quasi-linear
Schr\"odinger equations (see~\cite{CJS1}).
\vskip2pt
The reasons indicated above motivate the introduction of Proposition~\ref{abstract-result}.

\noindent
In the following Section~\ref{mainstatm} we present the statements of
the applications of the method indicated by this proposition
to four classes of constrained semi-linear and quasi-linear elliptic problems (more precisely, see
subsections~\ref{App1},~\ref{App2},~\ref{App3} and~\ref{App4}). Finally, in Section~\ref{proofsect}
we provide the proofs of the results stated in Section~\ref{mainstatm} (see, respectively, the
subsections~\ref{section2},~\ref{section3},~\ref{section4} and~\ref{section5}).

\medskip

{\bf Acknowledgements:} The first author would like to thank A.
Farina and B. Sirakov for stimulating discussions. The authors
also thank H. Hajaiej for some useful comments on the paper.

\vskip18pt
\begin{center}\textbf{Notations.}\end{center}
\begin{enumerate}
\item For $N\geq 1$, we denote by $|\cdot|$ the euclidean norm in $\R^N$.
\item $\R^+$ (resp.\ $\R^-$) is the set of positive (resp.\ negative) real values.
\item For $p>1$ we denote by $L^p(\R^N)$ the space of measurable functions $u$ such that
$\int_{\R^N}|u|^pdx<\infty$. The norm $(\int_{\R^N}|u|^pdx)^{1/p}$ in $L^p(\R^N)$ is denoted by $\|\cdot\|_p$.
\item We denote by $L^\infty(\R^N)$ the set of bounded measurable functions endowed with the supremum
norm $\|\cdot\|_{\infty}=\sup_{x\in\R^N}|u(x)|$.
\item For $s\in\N$, we denote by $H^s(\R^N)$ the Sobolev space of functions $u$ in $L^2(\R^N)$
having generalized partial derivatives $\partial_i^k u$ in $L^2(\R^N)$ for all $i=1,\dots, N$ and
any $0\leq k \leq s$.
\item The norm $(\int_{\R^N}|u|^2dx+ \int_{\R^N}|\nabla u|^2dx)^{1/2}$ in $H^1(\R^N)$ is denoted by $\|\cdot\|$ and more generally, the norm in $H^s$ is denoted by $\| \cdot\|_{H^s}$.
\item We denote by $C_0^{\infty}(\R^N)$ the set of smooth and compactly supported functions in $\R^N$.
\item We denote by $B(x_0,R)$ a ball in $\R^N$ of center $x_0$ and radius $R$.
\end{enumerate}
\medskip

\section{Statements of the main results}
\label{mainstatm}

In this section we shall exhibit four examples in which we can successfully apply the
approach of Proposition~\ref{abstract-result} to constrained semi-linear and quasi-linear problems.

\subsection{A Choquard type problem in $\R^3$}
\label{App1}

We consider a variant of the classical Choquard
Problem (cf.\ \cite{lieb,PLchoq}). Precisely,
we minimize the functional $J:H\to\R$ defined by
\begin{equation}
    \label{choquard0}
J(u)=\int_{\R^3} j(u,|\nabla
u|)dx-\iint_{\R^6}\frac{u^2(x)u^2(y)}{|x-y|}dxdy\,\, \quad \text{over $\|u\|_{L^2(\R^3)}^2=c$},
\end{equation}
where $c$ is a fixed positive number. Here $H$ is given by $H^1(\R^3)$, and we assume that
$$
j:\R\times [0,\infty)\to\R^+,
$$
is continuous, convex and increasing with respect to the second
argument and that there exists $\nu>0$ such that
\begin{equation}
    \label{coerc1}
j(s,|\xi|)\geq\nu|\xi|^2,\quad \text{for all $s\in \R^+$ and all
$\xi\in\R^3$.}
\end{equation}
Moreover, there exists a positive constant $C$ such that
\begin{equation}
    \label{growthinfinity}
j(s,|\xi|) \leq C|s|^6 + C|\xi|^2,\quad \text{for all $s\in \R^+$
and all $\xi\in\R^3$.}
\end{equation}
Finally, we assume that
\begin{equation}
    \label{opposite1}
j(-s,|\xi|)\leq j(s,|\xi|), \quad\text{for all $s\in \R^-$ and
all $\xi\in\R^3$}.
\end{equation}
For all $c>0$, let us set
$$
m(c)=\min_{\|u\|_{L^2(\R^3)}^2=c} J(u).
$$
Our result is the following
\begin{proposition}\label{choquard}
    Under the assumptions~\eqref{coerc1}-\eqref{opposite1}, $m(c)$ is reached
    for all $c>0$.
\end{proposition}

Here the functional~\eqref{choquard0} is invariant under translations in $\R^3$
and, thus, the problem at infinity coincides with the given problem.
If one wants to treat this minimization problem using directly the
Compactness Concentration Principle of~\cite{lions1} one faces the
problem of checking the strict inequalities~\eqref{strict-inequalities}.  To
achieve this, one usually establish (see Lemma
II.1 of~\cite{lions1}) that
\begin{equation}\label{surlinear}
m(\theta \lambda) < \theta m(\lambda), \quad \text{for all
$\lambda \in ]0,c[$ and $\theta \in ]1, c/\lambda]$}.
\end{equation}
Under our assumptions on the Lagrangian $j(s, |\xi|)$ there is no reasons
for inequality~\eqref{surlinear} to be true. However we shall prove that
(H0)-(H1) hold and since $m_{\infty}(\lambda) = m(\lambda)<0$ for
any $\lambda \in ]0,c],$ that also condition (H2) is true. In order to check
(H1) we choose a minimizing sequence consisting of Schwarz
symmetric functions. The possibility to take a minimizing sequence
of this type, for general $j(s,|\xi|)$, has recently been
established in~\cite{HaSq} for even weaker growth assumptions on $j$.
\smallskip

\subsection{A general class of quasi-linear problems}
\label{App2}
We study a general problem of minimization that goes back to the work of
Stuart~\cite{stuarbbelow} and has recently undergone new
developments~\cite{HaSq}. Let
\begin{equation}
    \label{minprob}
T=\inf\big\{J(u):\, u\in {\mathcal C}\big\},
\end{equation}
where we have set
\begin{equation*}
    {\mathcal C}=\Big\{u\in H:\text{$G_k(u_k),\,
j_k(u_k,|\nabla u_k|)\in L^1(\R^N)$ for any $k$ and $\sum_{k=1}^m
\int_{\R^N}G_k(u_k)dx=1$}\Big\},
\end{equation*}
being $m\geq 1$ and $H=W^{1,p}(\R^N,\R^m)$.
Here $J$ is a functional defined, for any function $u=(u_1,\dots,u_m)\in {\mathcal C}$, by
$$
J(u)=\sum_{k=1}^m \int_{\R^N} j_k(u_k,|\nabla u_k|)dx-\int_{\R^N}
F(|x|,u_1,\dots,u_m)dx.
$$
We collect below the assumptions on $j_k,F,G$
that we shall need to state the result.
\vskip3pt
\noindent
$\bullet$ {\bf Assumptions on $j_k$.}
For $m\geq 1$, $N\geq 1$, $p>1$, let
$$
j_k:\R\times [0,\infty)\to\R^+,\quad \text{for $k=1,\dots,m$}
$$
be continuous, convex and increasing functions with respect to the
second argument and such that there exists $\nu>0$ with, for
$k=1,\dots,m$,
\begin{equation}
    \label{coerc}
\nu|\xi|^p\leq j_k(s,|\xi|),\quad \text{for all $s\in \R^+$ and
all $\xi\in\R^N$.}
\end{equation}
Moreover there exist $\alpha >0$ and $\beta >0$ such that
\begin{equation}
    \label{coerc2}
 j_k(s,|\xi|) \leq \beta|\xi|^p,\quad \text{for all $s\in
[0, \alpha]$ and all $\xi\in\R^N$ with $|\xi| \in [0, \alpha]$.}
\end{equation}
Finally we require, for $k=1,\dots,m$,
\begin{equation}
    \label{opposite}
\text{$j_k(-s,|\xi|)\leq j_k(s,|\xi|)$,
\,\, for all $s\in \R^-$ and all $\xi\in\R^N$}.
\end{equation}

\vskip3pt
\noindent
$\bullet$ {\bf Assumptions on $F$.}
Let us consider a function
$$
F:[0,\infty)\times\R^m\to\R,
$$
of variables $(r,s_1,\dots,s_m)$, measurable and bounded with respect $r$ and
continuous with respect to $(s_1,\dots,s_m)\in\R^N$ with
$F(r,0,\dots,0)=0$ for any $r \in \R^+$. We assume that
\begin{align}
& F(r, s + he_i + ke_j)+ F(r, s) \geq F(r, s + he_i)+ F(r, s + ke_j), \label{supermod1}\\
& F(r_1,s + he_i)+ F(r_0,s) \leq F(r_1,s)+ F(r_0,s + he_i),   \label{supermod2}
\end{align}
for every $i\neq j$, $i,j=1,\dots,m$ where $e_i$ denotes the
$i$-th standard basis vector in $\R^m$, $r > 0$, for all $h,k >
0$, $s =(s_1,\dots,s_m)$ and $r_0,r_1$ such that $0 <r_0 <r_1$.
\vskip2pt Conditions~\eqref{supermod1}-\eqref{supermod2} are also
known as cooperativity conditions. Also, if $F$ is
smooth,~\eqref{supermod1} yields $\partial^2_{ij}F(r,s_1,\dots,s_m)\geq 0$ for $i\neq j$.
In general,~\eqref{supermod1}-\eqref{supermod2} are necessary for rearrangement inequalities to hold
(see~\cite{troy}). Moreover, we assume that
\begin{align}
&\limsup_{(s_1,\dots,s_m)\to (0,\dots,0)^+} \frac{F(r,s_1,\dots,s_m)}{\sum\limits_{k=1}^m s_k^{p}}<\infty,\label{zerocv} \\
&\lim_{|(s_1,\dots,s_m)|\to \infty}
\frac{F(r,s_1,\dots,s_m)}{\sum\limits_{k=1}^m
s_k^{p+\frac{p^2}{N}}}=0, \label{grothF}
\end{align}
uniformly with respect to $r$.

For a $j \in \{1,  \dots, m\}$ there exist $r_0>0$, $\delta>0$,
$\mu>0$, $\tau\in[0,p)$ and $\sigma\in[0,\frac{p(p-\tau)}{N}[$
such that $F(r,s_1,\dots,s_m) \geq 0$ for $|s| \leq \delta$ and
\begin{equation}
F(r,s_1,\dots,s_m)\geq  \mu r^{-\tau}s_j^{\sigma+p},
\qquad\text{for $r>r_0$ and $s\in\R^m_+$ with $|s|\leq\delta$}.
\label{zerocvBis}
\end{equation}
Also,
\begin{equation}
 \lim_{\underset{(s_1,\dots,s_m)\to (0,\dots,0)^+}{r\to+\infty}}
\frac{F(r,s_1,\dots,s_m)}{\sum\limits_{k=1}^m s_k^{p}}=0.
\label{zeroultima}
\end{equation}
Finally, we require:
\begin{equation}
    \label{modineq}
    F(r,s_1,\dots,s_m)\leq F(r,|s_1|,\dots,|s_m|), \quad \text{for all $r>0$ and
$(s_1,\dots,s_m)\in \R^m$}
\end{equation}
and for a $j \in \{1,\dots,m\}$ and a $\delta >0$
\begin{equation}\label{peccato}
s_{j} \to F(r,s_1,\dots, s_j,\dots, s_m) \quad \text{ is strictly
increasing for $s_j \in [0, \delta]$}.
\end{equation}

\vskip3pt
\noindent
$\bullet$ {\bf Assumptions on $G_k$.}
Consider $m\geq 1$ continuous  functions
$$
G_k:\R\to\R^+,\quad G_k(0)=0,\quad \text{for $k=1,\dots,m$}
$$
such that there exists $\gamma>0$ with
\begin{equation}
    \label{Gkass}
G_k(s)\geq \gamma |s|^p,\quad \text{for all $s\in\R$}.
\end{equation}
We also require
\begin{equation}
    \label{Gkass2}
G_j \quad \text{is $p$-homogeneous where $j \in \{1,  \dots, m\}$
is defined in \eqref{zerocvBis}}.
\end{equation}
\vskip2pt
Under the assumptions~\eqref{coerc}-\eqref{Gkass2}, we prove the following
\begin{theorem}
    \label{mainappl}
    Assume that $N=1$ and that~\eqref{coerc}-\eqref{Gkass2} hold.
Then problem~\eqref{minprob} admits a radially symmetric and
radially decreasing nonnegative solution. Furthermore for $N\geq 1$, if~\eqref{zerocvBis}
holds with $\tau=0$ and~\eqref{coerc2} holds for
all $s\in\R^+$ and $\xi\in\R^N$, then the same conclusion holds without
condition~\eqref{peccato}.
\end{theorem}

In problem~\eqref{minprob}, (H0) naturally hold and also (H1)
since we can choose a suitable minimizing sequence consisting of Schwarz
symmetric functions as in Section~\ref{section3}. Our effort here
is thus to derive weak assumptions under which condition (H2) is fulfilled. Taking
advantage that the minimizing sequence consists of radially
symmetric functions we can check (H2) constructing explicitly a mass $v$
such that $u+v \in {\mathcal C}$ and $J(u+v) <J(u)$.
In the first part of the statement, we restrict to $N=1$ since in checking (H2)
we use geometric properties of the graph of elements of $H^1(\R)$.
It is an open question if our result also holds for $N \geq 2$
(see also Proposition~\ref{restrictions} in Section~\ref{App3}).

\begin{remark}
    \label{no-growth}
In~\cite{HaSq} (see also~\cite{stuarbbelow}), in order to prove
that the weak limit $u$ satisfies the constraint, the growth of $j_k$
is related to the one of $F(|x|,s_1,\cdots,s_m)$. More precisely,
in~\cite{HaSq} it is assumed that there exists $\alpha\geq p$ such that
\begin{equation}
   \label{antimonot}
\text{$j_k(ts,t|\xi|)\leq t^\alpha j_k(s,|\xi|)$, \,\,\, for all
$t\geq 1$, $s\in \R^+$ and $\xi\in\R^N$}.
\end{equation}
and
\begin{equation}
    \label{psphere}
    F(r,t s_1,\dots,t s_m)\geq t^{\alpha} F(r,s_1,\dots,s_m),
\end{equation}
for all $r>0$, $t\geq 1$ and $(s_1,\dots,s_m)\in \R^m$, where
$\alpha\geq p$ is the value appearing in
condition~\eqref{antimonot}. Note that under~\eqref{antimonot}
and~\eqref{psphere} one has
$$
m(\lambda c) \leq \lambda^{\alpha}m(c),\quad \text{for any $c>0$ and $\lambda \geq 1$}.
$$
In particular $c \mapsto m(c)$ is strictly decreasing and
Proposition~\ref{Proposition2} yields the assertion.
\end{remark}

\begin{remark}
Take $\beta\geq 0$, $\tau\in[0,p)$, $\sigma\in[0,\textstyle{\frac{p(p-\tau)}{N}}]$
and a continuous and decreasing function $a:[0,\infty)\to[0,\infty)$ such that
$$
a(|x|)={\mathcal O}\left(|x|^{-\tau}\right)\quad\text{as $|x|\to\infty$}.
$$
Then the function
    $$
    F(|x|,s_1,\dots,s_m)=\frac{a(|x|)}{p+\sigma}\sum_{k=1}^m |s_k|^{p+\sigma}+
    \frac{2\beta a(|x|)}{p+\sigma}
    \sum_{\overset{i,j=1}{i\neq j}}^m|s_i|^{\frac{p+\sigma}{2}}|s_j|^{\frac{p+\sigma}{2}}
    $$
    satisfies all the required assumptions.
\end{remark}
\medskip

\subsection{A Stuart's type problem}
\label{App3}
We consider here the problem
\begin{eqnarray} \label{0.11}
\mbox{ minimize} \quad I \quad \mbox{ on } \quad
\|u\|_{L^2}^2= c
\end{eqnarray}
where  $c>0$ and $I : H^1(\R^N) \to \R$ is given by
$$
I(u) = \frac{1}{2}\int_{\R^N} |\nabla u|^2 dx  - \int_{\R^N}
F(x,u)dx.
$$
We discuss problem~\eqref{0.11} under the assumptions:
\begin{equation}\label{zerocvs}
\limsup_{s \to 0^+} \frac{F(x,s)}{s^2}<\infty  \quad \mbox{ and }
\quad \lim_{s\to \infty} \frac{F(x,s)}{s^{2+\frac{4}{N}}}=0,
\end{equation}
uniformly with respect to $x \in \R^N$. Also
\begin{equation}\label{compact}
\lim_{|x|\to \infty} F(x,s) =0,
\quad\text{uniformly in $s \in \R$},
\end{equation}
\begin{equation}
    \label{modineqs}
    F(x,s)\leq F(x,|s|),
\quad \text{for all $x \in \R^N$ and $s\in \R$.}
\end{equation}

\begin{remark}
    In some cases, for instance when $F$ has the form $F(x,s)=r(x)G(s)$ for any $x\in\R^N$ and $s\in\R$,
    assumption~\eqref{compact} can be relaxed by just asking that $r(x)\to 0$ as $|x|\to\infty$.
\end{remark}

In addition, we consider the following assumption: there exists a
positive constant $\delta$ such that $F: \R^N \times [0, \delta]
\to \R^+$ is a Carath\'{e}odory function and
\begin{equation}
    \label{highdimension}
\begin{cases}
\text{$N \geq 1$  and there exist  $r_0,A>0$, $d \in (0,2)$
   and $\alpha \in (0,\frac{2(2-d)}{N})$ with}  \\
 F(x,s) \geq A (1 +|x|)^{-d} s^{2 + \alpha},\quad \mbox{ for all } s
\in [0, \delta] \mbox{ and   } |x| \geq r_0, & \\
\noalign{\vskip4pt} \text{$N=1$ and there exist $r_0 >0$ and
$\alpha \in (0,2)$
with} \\
 F(x,s) \geq r(x) s^{2 + \alpha},\quad  \mbox{ for all } s \in
[0, \delta] \mbox{ and   } |x| \geq r_0, &
\end{cases}
\end{equation}
where $r \in L^{\infty}(\R), \ r\geq 0$  and
$$
\int_{\R\backslash [-r_0,r_0]}r(x)dx >0,
$$
where the value $+ \infty$ is admissible.

\begin{remark}
If we consider problem~\eqref{0.11} within the formalism
of~\cite{lions1} we see that, because of~\eqref{compact}  the
associated ``problem at infinity" is
\begin{equation*}
\text{minimize \,\,\, $I_{\infty}(u)= \frac{1}{2}\int_{\R^N}
|\nabla u|^2 dx$\,\,\,\, on $\|u\|_2^2= c$}.
\end{equation*}
Thus setting
$$
m_{\infty}(c) =  \inf  \{ I_{\infty}(u): \text{$u\in H^1(\R^N)$ with $\|u\|_{2}^2 = c$}\},
$$
we have $m_{\infty}(c)=0$.
\end{remark}

Assumptions~\eqref{zerocvs}-\eqref{highdimension} are classical
assumptions first introduced in~\cite{stuarbbelow} under which $I$
is well defined and continuous. Also (H0) is known to hold and,
because of~\eqref{compact}, any minimizing sequence
for~\eqref{0.11} satisfies (H1). Now defining
$$
m(c) =  \inf  \{ I(u):\, \|u\|_{2}^2 = c\},
$$
we have the following

\begin{proposition}
    \label{qualitative}
Assume that~\eqref{zerocvs}-\eqref{highdimension} hold. Then $m(c)
< 0$ for all $c > 0$ and $c \mapsto m(c)$ is non increasing.
\end{proposition}

\begin{remark}\label{crucial}
Assume that conditions~\eqref{zerocvs}-\eqref{highdimension} hold
and let $u\in H^1(\R^N)$ be a function such that $\|u\|_2^2\leq c$
and $I(u) \leq m(c) <0$ (such a $u$ comes from a weakly convergent
minimizing sequence $(u_n)$ over which the functional $I$ is lower
semicontinuous). Then $u \in H^1(\R^N)$ minimizes $I$ on the
constraint $d:= \|u\|_2^2>0$.  Indeed if there exists $v \in
H^1(\R^N)$ with $\|v\|_2^2 = \|u\|_2^2=d$ and $I(v) < I(u)$ we get
a contradiction since, by Proposition~\ref{qualitative}, the map
$\lambda \mapsto m(\lambda)$ is non increasing.
\end{remark}

To show that $m(c)$ is reached we must restrict our assumptions.
First we have

\begin{proposition}\label{restrictions}
Assume that ~\eqref{zerocvs}-\eqref{highdimension}  hold. In
addition assume that $N=1$ and there exists $\delta >0$ such that,
for any $x\in\R$,
\begin{equation}
    \label{increasing2}
\text{$s \mapsto F(x,s)$ is strictly increasing for $s \in [0,\delta]$}.
\end{equation}
Then $m(c)$ is reached.
\end{proposition}

Our second result requires some additional regularity of the
nonlinearity $F(x,s)$. We assume that the derivative $f(x,s)=
F_s(x,s)$  of $F(x,s)$ with respect to $s \in \R$ exists, that
$f:\R^N\times\R^+\to\R^+$ is a Carath\'{e}odory function and satisfy
\begin{equation}\label{zerocvsbis}
\limsup_{s \to 0^+} \frac{f(x,s)}{s}<\infty  \quad \mbox{ and }
\quad \lim_{s\to +\infty} \frac{f(x,s)}{s^{1+\frac{4}{N}}}=0,
\end{equation}
uniformly with respect to $x \in \R^N$. We also
replace~\eqref{highdimension} by
\begin{equation}
    \label{highdimensionbis}
\begin{cases}
\text{$N <5 $  and there exist  $r_0,A>0$, $d \in (0,2)$
   and $\alpha \in (0,\frac{2(2-d)}{N})$ with}  \\
 f(x,s) \geq A (1 +|x|)^{-d} s^{1 + \alpha},\quad \mbox{ for all } s
\in \R^+ \mbox{ and   } |x| \geq r_0, & \\
\noalign{\vskip4pt} \text{$N\geq 5$ and there exist $r_0,A>0$, $d \in (0,2)$
   and $\alpha \in (0,\frac{2-d}{N-2})$ with}  \\
 f(x,s) \geq A (1 +|x|)^{-d} s^{1 + \alpha},\quad \mbox{ for all } s
\in \R^+ \mbox{ and   } |x| \geq r_0. &
\end{cases}
\end{equation}

\begin{proposition}\label{restrictionsbis}
Assume that  \eqref{compact}-\eqref{modineqs}
and~\eqref{zerocvsbis}-\eqref{highdimensionbis} hold. Then $m(c)$ is reached.
\end{proposition}

\subsection{A problem studied by Badiale-Rolando}
\label{App4}
Finally, we consider in this section the following
problem: Let $x=(y,z) \in \R^k \times \R^{N-k}$ with $N >k \geq 2$
and set
\begin{align*}
&H:=\Big\{ u \in H^1(\R^N): \,\, \int_{\R^N} \frac{|u|^2}{|y|^2}dx <
\infty\Big\}  \\
& H_s := \Big\{ u \in H: \,\, u(y,z) =
u(|y|,z)\Big\}.
\end{align*}
Let $f: \R \to \R$ be continuous and satisfies, for $F(t) :=
\int_0^t f(s) ds,$
\begin{itemize}
\item[$(f_0)$] $F(t_0) >0$ for some $t_0 >0$.
\item[($f_1)$] there exists $q>2$ such that
$$
\lim_{t\to 0^+}\frac{f(t)}{|t|^{q-1}}=0,
$$
\end{itemize}
and one of the following assumptions:
\begin{itemize}
\item[$(f_2)$] $f(\beta)=0$ for some $\beta > \beta_0 := \inf\{t >0,
F(t) >0\}.$
\item[($f_3)$] there exists $p \in ]2, 2 + \frac{4}{N}[$ such that
$$
\lim_{t\to+\infty} \frac{f(t)}{|t|^{p-1}}=0.
$$
\end{itemize}
\vskip2pt
\noindent
Our result is stated in the following
\begin{theorem}
    \label{Badiale-Rolando}
    Let $N > k \geq 2$ and $\mu >0$. Assume that $f \in C(\R, \R)$
    satisfies $(f_0), (f_1)$ and at least one of the hypotheses
    $(f_2)$ and $(f_3)$. Then there exists $\rho_0 >0$ such that
    for all $\rho > \rho_0$ the minimization problem
    \begin{equation}
    \label{minimization-problem}
    \inf_{u \in H_s, \, \|u\|_2^2 = \rho} \Big(
    \frac{1}{2}\int_{\R^N} |\nabla u|^2 dx + \frac{\mu}{2}
    \int_{\R^N}\frac{|u|^2}{|y|^2}dx - \int_{\R^N} F(u) dx \Big)
    \end{equation}
admits a solution $u(y,z) = u(|y|,|z|) \geq 0$ which is non
increasing in $|z|$.
\end{theorem}

Theorem~\ref{Badiale-Rolando} was originally proved
in~\cite{BaRo}. It is the central part of~\cite{BaRo} in which is
establish the existence of standing waves with non zero angular
momentum for a class of Klein-Gordon equations. We refer
to~\cite{BaRo} for a detailed presentation of the problem and of
its physical motivations. Here we concentrate on giving an
alternative shorter proof of this result. The original proof
in~\cite{BaRo} is based on the full machinery of the Concentration
Compactness Principle and the central issue is to rule out the
dichotomy case. Here we follow the added-mass approach presented in
Proposition~\ref{abstract-result}. Due to the symmetry
of~\eqref{minimization-problem} it is possible to choose a
minimizing sequence such that (H1) holds. Then, still using the
symmetry, a simple scaling argument shows that (H2) holds as well.
\bigskip

\section{Proofs of the main results}
\label{proofsect}
In the following section we prove all the achievements announced in Section~\ref{mainstatm}.

\subsection{Proof of Proposition~\ref{choquard}}
\label{section2}
We define the Coulomb energy in $\R^3$ by setting
$$
\D(u)=\iint_{\R^6}\frac{u^2(x)u^2(y)}{|x-y|}dxdy,
$$
for all $u\in H^1(\R^3)$. First we have the following

\begin{lemma}\label{continuity-choquard}
\label{finitEn} Let $u \in H^1(\R^3)$ with
$\|u\|_{L^2(\R^3)}^2=c>0$. There exists a positive constant $C$,
depending only on $c$, such that
\begin{equation*}
\D(u) \leq C \|u\|_{H^1(\R^3)}.
\end{equation*}
\end{lemma}


\begin{proof}
Combining Hardy-Littlewood-Sobolev inequality (see e.g.\
Lieb-Loss, Thm 4.3, p.106) with Gagliardo-Nirenberg inequality,
yields a positive constant $C_0$ such that
\begin{equation}
    \label{Hardyineq}
\D(u) \le C_0\|u\|^4_{L^{\frac{12}{5}}(\R^3)} \le
C_0\|u\|^{3}_{L^2(\R^3)}\,\|u\|_{H^1(\R^3)} =C_0
c^{3/2}\,\|u\|_{H^1(\R^3)},
\end{equation}
which concludes the proof.
\end{proof}

Secondly, we need the following approximation result.

\begin{lemma}\label{compact-support}
 Assume that conditions~\eqref{coerc1}-\eqref{opposite1} hold. Let $u \in
 H^1(\R^3) \backslash\{0\}$ be given. Then, for any $\varepsilon >0$
 there exists $\tilde{u} \in C_0^{\infty}(\R^3)$ such that
 $$
J(\tilde{u}) \leq J(u) + \varepsilon \quad \mbox{and} \quad
 \|\tilde{u}\|_{L^2(\R^3)}^2 = \|u\|_{L^2(\R^3)}^2.
$$
 \end{lemma}

\begin{proof}
By density of $C_0^{\infty}(\R^3)$ into $H^1(\R^3)$ there exists a sequence
$(u_n) \subset C_0^{\infty}(\R^3)$ with $u_n \to u$ in
$H^1(\R^3)$, as $n\to\infty$. In particular
$\|u\|_{L^2(\R^3)}/\|u_n\|_{L^2(\R^3)} \to 1$, as $n\to\infty$. Thus
$$
\left\|u- \frac{\|u\|_{L^2(\R^3)}}{\|u_n\|_{L^2(\R^3)}}u_n\right\|_{H^1(\R^3)}
\leq \|u-u_n\|_{H^1(\R^3)} + \left|1 -
\frac{\|u\|_{L^2(\R^3)}}{\|u_n\|_{L^2(\R^3)}}\right| \|u_n\|_{H^1(\R^3)} \to 0,
$$
as $n\to\infty$. This
proves that there exists a sequence $(\tilde{u_n}) \subset
C_0^{\infty}(\R^3)$ with $\|\tilde{u_n}\|_{L^2(\R^3)} = \|u\|_{L^2(\R^3)}^2$ such that
$\tilde{u_n} \to u$ in $H^1(\R^3)$. To conclude we just need to
prove that $J(\tilde{u_n}) \to J(u)$, as $n\to\infty$. Clearly, by
Lemma~\ref{continuity-choquard}, $\D(\tilde{u_n})\to \D(u)$ (see e.g.\ estimate~\eqref{stimaDD} hereafter). Now
from the growth condition~\eqref{growthinfinity}, by the
generalized Lebesgue Theorem (see Theorem IV of~\cite{c2}) we
readily get that  $\int_{\R^3}j(\tilde{u_n}, |\nabla
\tilde{u_n}|) dx \to \int_{\R^3}j(u, |\nabla u|) dx$, as $n\to\infty$.
\end{proof}

We can now give the proof of Proposition~\ref{choquard}.

\begin{proof}
    Let us fix a positive number $c$ and let $(u_h)\subset H^1(\R^3)$ be a minimizing sequence for $m(c)$,
    namely $\|u_h\|^2_2=c$, for all $h\geq 1$, and
    $$
    \int_{\R^3} j(u_h,|\nabla u_h|)dx=m(c)+\D(u_h)+o(1),\quad
    \text{as $h\to\infty$}.
    $$
    By virtue of Lemma~\ref{finitEn} and assumption~\eqref{coerc1}, we have
    $$
    \nu\|\nabla u_h\|_{L^2(\R^3)}^2\leq m(c)+C\|\nabla u_h\|_{L^2(\R^3)}+o(1),\quad
    \text{as $h\to\infty$},
    $$
    so that $(u_h)$ is bounded in $H^1(\R^3)$ and assumption (H0)
    of  Proposition~\ref{abstract-result} is thus satisfied. Up to a subsequence, $(u_h)$ weakly
    converges to some function $u$ in $H^1(\R^3)$.
    Observe now that, if $u^*_h$ denotes the symmetrically decreasing
    rearrangement of $u_h$, for all $h\geq 1$,
    $$
    \iint_{\R^6}\frac{u^2_h(x)u^2_h(y)}{|x-y|}dxdy\leq\iint_{\R^6}\frac{(u^2_h)^*(x)(u^2_h)^*(y)}{|x-y|}dxdy=
    \iint_{\R^6}\frac{(u^*_h)^2(x)(u^*_h)^2(y)}{|x-y|}dxdy,
    $$
    where we have used the fact that $(u^*_h)^2=(u^2_h)^*$.
    For this rearrangement inequality, started with the work of Lieb~\cite{lieb}, we refer
for instance to~\cite{burchhaj}.

In turn, by taking into account that by~\cite[Corollary 3.3]{HaSq} we have
$$
\int_{\R^3} j(u^{*}_h,|\nabla u^{*}_h|)dx\leq \int_{\R^3}
j(u_h,|\nabla u_h|)dx,
$$
we conclude that $J(u_h^*)\leq J(u_h)$, for all $h\geq 1$. Hence,
we may assume that $(u^*_h)$ is a positive (since $J(|v|)\leq
J(v)$, for all $v\in H^1(\R^3)$) minimizing sequence for $J$,
which is radially symmetric and radially decreasing. In what
follows, we denote it again by $(u_h)$. Taking into account that
$(u_h)$ is bounded in $L^{2}(\R^3)$, it follows that (see
\cite[Lemma A.IV]{BL1}) $u_h(x)\leq M|x|^{-3/2}$ for all
$x\in\R^3\setminus\{0\}$ and $h\in\N$, for some constant $M>0$ and hence
$(u_h)$ turns out to be strongly convergent to $u$
in $L^q(\R^3)$ for any $2<q<6$. In particular, we have the strong limit
\begin{equation}
    \label{conver125}
u_h\to u\quad\text{in $L^{\frac{12}{5}}(\R^3)$,\, as
$h\to\infty$}.
\end{equation}
We  want to show that
$$
\D(u_h)\to \D(u),\quad\text{as $h\to\infty$}.
$$
To this end, we use that the Coulomb potential $|x|^{-1}$ is even
and write
\begin{equation*}
|\D(u_h)-\D(u)|\leq \D(||u_h|^2
-|u|^2|^{1/2},(|u_h|^2+|u|^2)^{1/2}).
\end{equation*}
Let us now introduce the two variable functional
$$
\D(v,w):=\iint_{\R^6}\frac{v^2(x)w^2(y)}{|x-y|}dxdy,
$$
for all $v,w\in H^1(\R^3)$. The following inequality holds (see
e.g.\ Lieb-Loss, Thm 9.8, p.250)
\begin{equation}
    \label{shwa}
\D(v,w)^2\le \D(v,v) \,\D(w,w),\quad\text{for all $v,w\in
H^1(\R^3)$}.
\end{equation}
Now, by means of Hardy-Littlewood-Sobolev inequality (see the
first line of~\eqref{Hardyineq}) as well as
H\"older's inequality,  it follows that (just use
inequality~\eqref{shwa} with $v=v_h=||u_h|^2-|u|^2|^{1/2}$ and
$w=w_h=(|u_h|^2+|u|^2)^{1/2}$ for all $h\geq 1$) there exists a
constant $C$ with
\begin{align}
    \label{stimaDD}
|\D(u_h)-\D(u)|^2 &\leq  C
\|\,||u_h|^2-|u|^2|^{1/2}\|_{L^\frac{12}{5}(\R^3)}^{4}\|
\,(|u_h|^2+|u|^2)^{1/2}\|_{L^\frac{12}{5}(\R^3)}^{4} \\
&\leq  C \|u_h-u\|_{L^\frac{12}{5}(\R^3)}^{2}.  \notag
\end{align}
This implies, via~\eqref{conver125}, the desired convergence of
$\D(u_h)$ to $\D(u)$. Also as $j(s,t)$ is positive, convex and
increasing in the second argument (and thus $\xi\mapsto
j(s,|\xi|)$ is convex), $u_h\to u$ in $L^1_{{\rm loc}}(\R^3)$ and
$\nabla u_h\rightharpoonup  \nabla u$ in $L^1_{{\rm loc}}(\R^3)$,
by well known lower semicontinuity results (cf.\
\cite{ioffe1,ioffe2}) it follows that
\begin{equation}\label{los}
 \int_{\R^N}j(u,|\nabla u|)dx\leq
\liminf_{h}\int_{\R^N}j(u_h,|\nabla u_h|)dx,
\end{equation}
and we can conclude that
$$
J(u)\leq\liminf_{h\to\infty} J(u_h).
$$
Therefore, also condition (H1) is fulfilled.

Now, given a function $w\in C_0^{\infty}(\R^3)$ with
$\|w\|_2^2=c$, and considering the rescaling $\{t\mapsto w_t\}$
with $w_t(x)=t^{3/2}w(tx)$, we have $\|w_t\|_2^2=c$ for all $t>0$
and
$$
\D(w_t)=\iint_{\R^6}\frac{w_t^2(x)w_t^2(y)}{|x-y|}dxdy=
t^6\iint_{\R^6}\frac{w^2(tx)w^2(ty)}{|x-y|}dxdy=t\D(w).
$$
Hence, taking into account the growth condition~\eqref{growthinfinity}, we conclude
\begin{align*}
m(c) & \leq \int_{\R^3} j(w_t,|\nabla w_t|)dx-\D(w_t) \\
& \leq C\int_{\R^3} |w_t|^6 dx+C\int_{\R^3} |\nabla w_t|^2 dx-t\D(w) \\
&= Ct^6\int_{\R^3} |w|^6 dx+Ct^2\int_{\R^3} |\nabla w|^2 dx-t\D(w)<0,
\end{align*}
for  $t>0$ sufficiently small. In turn, we have $J(u) \leq
m(c)<0$, which also yields $u\neq 0$. Now, if it was $\|u\|_{L^2(\R^3)}^2 =c$, the
proof would be over. Otherwise we assume, by contradiction, that
$\|u\|_{L^2(\R^3)}^2=\lambda$ with $0<\lambda<c$. Following the proof that
$m(c) <0$, we see that there exists a function $v \in C_0^{\infty}(\R^3)$
such that $\|v\|_{L^2(\R^3)}^2 = c - \lambda >0$ and $J(v) <0$. Also by
Lemma~\ref{compact-support}, it is
possible to find a $\tilde{u} \in C_0^{\infty}(\R^3)$ with
$\|\tilde{u}\|_{L^2(\R^3)}^2 = \lambda$ and $J(\tilde{u}) \leq J(u) +
\frac{|J(v)|}{2}$. Taking advantage that~\eqref{choquard0} is an
autonomous problem we can assume that $v$ and $\tilde{u}$ have
disjoint supports. Thus
$$
\|v+ \tilde{u}\|_{L^2(\R^3)}^2 = \|v\|_{L^2(\R^3)}^2 + \|\tilde{u}\|_{L^2(\R^3)}^2 = (c-
\lambda) + \lambda = c,
$$
as well as
$$
J(v + \tilde{u}) = J(v) +
J(\tilde{u}) \leq J(v) + J(u) - \frac{J(v)}{2} \leq J(u) +
\frac{J(v)}{2} < J(u).
$$
Thus (H2) hold and the proof is completed.
\end{proof}
\smallskip

\subsection{Proof of Theorem~\ref{mainappl}}
\label{section3}
We shall divide the proof into three main steps.
The first part of the proof (Step I), aiming to prove that conditions (H0) and (H1)
of our abstract machinery hold, follows the pattern of
the proof of~\cite[Theorem 4.5]{HaSq}. For the sake of completeness
we report here some of the arguments in order to have a complete
picture of the situation. Instead, the last part of the
proof (Steps II and III) contains the main elements of novelty and
improvement (through to the mass addiction argument) with
respect to~\cite[Theorem 4.5]{HaSq}.
\vskip5pt
\noindent
{\bf Step I. [Verification of (H0) and (H1)]}
Let $u^h=(u^h_1,\dots,u^h_m)\subset {\mathcal C}$
be a minimizing sequence for the functional $J$. Then
\begin{gather}
 \lim_{h}\Big(\sum_{k=1}^m \int_{\R^N} j_k(u_k^h,|\nabla u_k^h|)dx-\int_{\R^N} F(|x|,u_1^h,\dots,u_m^h)dx\Big)=T, \label{min1}\\
 G_k(u^h_k),\, j_k(u_k^h,|\nabla u_k^h|)\in L^1(\R^N), \quad
\sum_{k=1}^m\int_{\R^N}G_k(u^h_k)dx=1,\quad\text{for all $h\in\N$}.\notag
\end{gather}
In light of~\eqref{opposite} and~\eqref{modineq}, we obtain
$J(|u^h_1|,\dots,|u^h_m|)\leq J(u^h_1,\dots,u^h_m)$ for all $h\in\N$,
so we may assume, without loss of generality, that $u^h_k\geq 0$
a.e.\ in $\R^N$, for all $k=1,\dots,m$ and $h\in\N$. Now one can prove that
$(u^h)$ is bounded in $W^{1,p}(\R^N,\R^m)$. To this aim, since
$(u^h)\subset {\mathcal C}$, by assumption~\eqref{Gkass} on $G_k$,
the sequence $(u^h)$ is bounded in $L^p(\R^N)$. By combining the
growths~\eqref{zerocv}-\eqref{grothF}, for every $\eps>0$
there exists $C_\eps>0$ with
\begin{equation}
    \label{subcrit}
F(r,s_1,\dots,s_m)\leq C_\eps\sum_{k=1}^m s_k^p+\eps \sum_{k=1}^m s_k^{p+\frac{p^2}{N}},\quad
\text{for all $r,s_1,\dots,s_m\in (0,\infty)$}.
\end{equation}
Therefore, in view of the Gagliardo-Nirenberg inequality
\begin{equation}
    \label{gagnir}
\|u^h_k\|_{L^{p+\frac{p^2}{N}}(\R^N)}^{p+\frac{p^2}{N}}\leq C\|u^h_k\|_{L^p(\R^N)}^{\frac{p^2}{N}}
\|\nabla u^h_k\|_{L^p(\R^N)}^{p},
\end{equation}
by combining~\eqref{coerc} with~\eqref{min1}, one immediately
yields the desired boundedness of $(u^h)$ in
$W^{1,p}(\R^N,\R^m)$. Hence condition (H0) hold for any positive
minimizing sequence.
\vskip3pt
\noindent
Now, after extracting a subsequence, still denoted by $(u^h)$, for any $k=1,\dots,m$,
\begin{equation}
\label{subseq}
u^h_k\rightharpoonup u_k\,\,\, \text{in $L^{p^*}(\R^N)$},\,\,
Du^h_k\rightharpoonup Du_k\,\,\, \text{in $L^{p}(\R^N)$},\,\,
u^h_k(x)\to u_k(x)\quad\text{a.e.\ $x\in\R^N$}.
\end{equation}
Of course, we have
$$
\sum_{k=1}^m\int_{\R^N}G_k(u_k)dx\leq\liminf_{h\to\infty}\sum_{k=1}^m\int_{\R^N}G_k(u^{h}_k)dx=1.
$$
In particular $G_k(u_k)\in L^1(\R^N)$.
For any $k=1,\dots,m$ and $h\in\N$, we denote by $u^{*h}_k$ the Schwarz symmetric rearrangement
of $u^h_k$. By means of~\cite[Theorem 1]{burchhaj}, we have
\begin{equation*}
\int_{\R^N} F(|x|,u_1^h,\dots,u_m^h)dx\leq \int_{\R^N} F(|x|,u_1^{*h},\dots,u_m^{*h})dx.
\end{equation*}
Moreover, by~\cite[Corollary 3.3]{HaSq}, we know that
$$
\int_{\R^N} j_k(u^{*h}_k,|\nabla u^{*h}_k|)dx\leq \int_{\R^N} j_k(u^h_k,|\nabla u^h_k|)dx.
$$
Finally, $u^{*h}\in{\mathcal C}$.
Hence, since $J(u^{*h})\leq J(u^h)$, $u^{*h}\in{\mathcal C}$, for $h\in\N$,
it follows that $u^{*h}=(u^{*h}_1,\dots,u^{*h}_m)$ is a positive
minimizing sequence for $J|_{{\mathcal C}}$, which is radially
symmetric and radially decreasing. In what follows, we denote it
again $u^{h}=(u^{h}_1,\dots,u^{h}_m)$. Taking into account that
$u^{h}_k$ is bounded in $L^{p}(\R^N)$, it follows that (see
\cite[Lemma A.IV]{BL1}) $u^h_k(x)\leq c_k|x|^{-N/p}$ for all
$x\in\R^N\setminus\{0\}$ and $h\in\N$,
for a positive constant $c_k$, independent of $h$. In turn,
by virtue of condition~\eqref{zeroultima}, for all
$\eps>0$ there exists $\rho_\eps>0$ such that
$$
|F(|x|,u_1^h(|x|),\dots,u_m^h(|x|))|\leq\eps\sum_{k=1}^m |u_k^h(|x|)|^{p},
\quad\text{for all $x\in\R^N$ with $|x|\geq \rho_\eps$}.
$$
Hence, it is easy to see that
\begin{equation*}
\int_{\R^N\setminus B(0,\rho_\eps)}F(|x|,u_1^h,\dots,u_m^h)dx
\leq\eps C,\quad
\int_{\R^N\setminus B(0,\rho_\eps)}F(|x|,u_1,\dots,u_m)dx\leq \eps C.
\end{equation*}
In turn, one readily obtains
\begin{equation}
\label{conFtoF}
\lim_h\int_{\R^N} F(|x|,u_1^h,\dots,u_m^h)dx =\int_{\R^N} F(|x|,u_1,\dots,u_m)dx.
\end{equation}
Also, arguing as in the proof of~\eqref{los}, for any
$k=1,\dots,m$ it follows
\begin{equation}
\label{lowersem}
\int_{\R^N}j_k(u_k,|Du_k|)dx\leq \liminf_{h}\int_{\R^N}j_k(u^h_k,|Du^h_k|)dx.
\end{equation}
Hence, $j_k(u_k,|Du_k|)\in L^1(\R^N)$ for any $k$ and
from~\eqref{conFtoF} and~\eqref{lowersem} it follows
\begin{equation}
    \label{lowerSC}
J(u)\leq \liminf_{h} J(u^h)=\lim_{h} J(u^h)=T.
\end{equation}
At this point also (H1) is established.

\vskip6pt
\noindent
{\bf Step II.}
To show that (H2) holds, let us first prove that $T<0$. For any
$\theta\in(0,1]$, we consider the function
$$
\Upsilon^\theta_j(x)=
\frac{\theta^{N/p^2}}{d_j^{1/p}}e^{-\theta|x|^p},\quad
d=\int_{\R^N} G_j(e^{-|x|^p})dx,
$$
where $j \in \{1,\dots,m\}$ is given by condition~\eqref{zerocvBis}. Without
restriction we can assume that $j=1$. Then by~\eqref{Gkass2} we
get
$$
\int_{\R^N}G_1(\Upsilon^\theta_1(x))dx=
\frac{\theta^{N/p}}{d}\int_{\R^N}G_1(e^{-\theta|x|^p})dx=
\frac{1}{d}\int_{\R^N}G_1(e^{-|x|^p})dx=1.
$$
Therefore $(\Upsilon^\theta_1,0,\dots,0)$ belongs to ${\mathcal
C}$ for any $\theta >0$.  Notice that
$$
|\nabla
\Upsilon^\theta_1(x)|^p=p^p\frac{\theta^{N/p+p}}{d}e^{-p\theta|x|^p}|x|^{p(p-1)}
,\quad x\in\R^N.
$$
Thus, for $\theta >0$ small enough, it follows  by~\eqref{coerc2},
that
\begin{align*}
\int_{\R^N} j_1(\Upsilon^\theta_1(x),|\nabla
\Upsilon^\theta_1(x)|)dx  &\leq \int_{\R^N} \beta |\nabla
\Upsilon^\theta_1(x)|^pdx \\
&\leq \frac{
p^p\theta^{N/p+p}}{d}\int_{\R^N}e^{-p\theta|x|^p}|x|^{p(p-1)}dx
=\theta C,
\end{align*}
where we have set
$$
C=\frac{ p^p}{d}\int_{\R^N}e^{-p|x|^p}|x|^{p(p-1)}dx.
$$
Now, in light of~\eqref{zerocvBis}, since taking $\theta >0$ small
enough we can assume that $0\leq \Upsilon^\theta_1\leq \delta$, we
obtain
$$ \int_{\R^N} F(|x|,\Upsilon^\theta_1(x),0,\dots,0)dx  \geq
\frac{\mu}{d^{\frac{\sigma +p}{p}}} \theta^{\frac{N(\sigma
+p)}{p^2}} \int_{\{|x|\geq\rho\}} |x|^{-\tau}
e^{-\theta(\sigma+p)|x|^p}dx \geq \theta^{\frac{N\sigma +p
\tau}{p^2}} C^{'}$$ with
$$
C^{'} = \frac{\mu}{d^{\frac{\sigma+p}{p}}}\int_{\{|x|\geq\rho\}}
|x|^{-\tau} e^{-(\sigma+p)|x|^p}dx.
$$
In conclusion, collecting the previous inequalities, for $\theta>0$ sufficiently small,
\begin{align*}
T &\leq  \int_{\R^N} j_1(\Upsilon^\theta_1(x),|\nabla
\Upsilon^\theta_1(x)|)dx
 -\int_{\R^N} F(|x|,\Upsilon^\theta_1(x),0, \dots,0)dx \\
& \leq \theta  \big(C-\theta^{\frac{N\sigma+p\tau-p^2}{p^2}}
C^{'}\big)<0,
\end{align*}
as $N\sigma+p\tau-p^2<0$, yielding the assertion. Notice that
$(u_1,\dots,u_m)\neq (0,\dots,0)$, otherwise we would get a
contradiction combining~\eqref{lowerSC} and $T<0$. We now define
$$
\zeta:=\sum_{k=1}^m\int_{\R^N}G_k(u_k)dx.
$$
If $\zeta=1$ then $( u_1,\dots, u_m)$ belongs to ${\mathcal C}$ and we
are done. We thus assume that $\zeta <1$ and look for a
contradiction.

\vskip6pt
\noindent
{\bf Step III-a. [Verification of (H2), $N\geq 2$, $\tau=0$ in~\eqref{zerocvBis}]}
Assuming that $\zeta <1$, we can conclude
as in Proposition~\ref{choquard}. Here~\eqref{minprob} is not
autonomous but the fact that~\eqref{zerocvBis} holds with $\tau=0$ permits to
select a $v\in C_0^{\infty}(\R^N)$ and a $\alpha>0$ such that
$$
\int_{\R^N}G_1(v) dx = 1 - \zeta, \quad
J(v(\cdot+y)) \leq - \alpha\quad \text{for any $y \in \R^N$ with  $|y|$ large enough}.
$$
Then we can conclude by arguing as in
Section~\ref{section2}, replacing the weak limit $u$ by
the compactly supported function $\tilde{u} \in C_0^{\infty}(\R^N)$  (cf.\ the proof of Lemma~\ref{compact-support}), thus avoiding the monotonicity condition~\eqref{peccato}.
\vskip6pt
\noindent
{\bf Step III-b. [Verification of (H2), $N=1$]}
Let $j \in \{1, \dots,m\}$ be such
that~\eqref{peccato} hold. Without restriction we can assume that
$j =1$. Since $u_1(x)$ is radially symmetric and positive we can
set $v_1(r) = u_1(|x|)$ with $v_1: \R^{+} \to \R^{+}$. We now
define $w_1 :\R \to \R^{+}$ by setting
$$
w_1(x):=
\begin{cases}
    v_1(|x|)   & \text{if $|x| \in [0,\varrho]$} \\
     v_1(\varrho)   & \text{if $|x| \in [\varrho,\varrho+ \mu]$}\\
    v_1(|x|-\mu)  & \text{if $|x| \in [\varrho+ \mu, \infty[$.}
\end{cases}
$$
Here $\varrho>0$ is such that $0 < v_1(\varrho) \leq \delta$ where $\delta
>0$ is given in condition~\eqref{peccato}. Without restriction we
can require $v_1$ to be continuous at $\varrho$. Instead, the value $\mu >0$ is
fixed in order to have
$$
\int_{[-\varrho-\mu,\varrho + \mu] \backslash [-\varrho,\varrho]} G_1(v_1(\varrho)) dx = 1-
\zeta.
$$
Now defining $w=(w_1,\dots,w_n) := (w_1,u_2,\dots,u_m)$ we
have by construction
\begin{equation} \nonumber
\sum_{k=1}^m\int_{\R}G_k( w_k)dx=1,
\end{equation}
namely $( w_1,\dots, w_m)$ belongs to the constraint ${\mathcal C}$.
Also, using~\eqref{coerc2},
\begin{equation} \nonumber
 \sum_{k=1}^m\int_{\R}j_k( w_k, |w'_k|)dx = \sum_{k=1}^m\int_{\R}j_k( u_k, |u'_k|)dx.
\end{equation}
Now split the integral as
$$
\int_{\R}F(|x|,w_1,\dots,w_n)dx =
\int_{[-\varrho,\varrho]}F(|x|,u_1,\dots,u_m)dx + \int_{\R^N \backslash
[-\varrho,\varrho]}F(|x|,w_1,\dots,w_m)dx.
$$
We have $w_1(x) \geq u_1(x)$ a.e.\ in $\R$ and $w_1
\neq u_1$, so recalling the monotonicity condition~\eqref{peccato} we have
\begin{equation} \nonumber
\int_{\R \backslash [-\varrho,\varrho]}F(|x|,w_1,\dots,w_m)dx >  \int_{\R
\backslash [-\varrho,\varrho]}F(|x|,u_1,\dots,u_m)dx.
\end{equation}
We then deduce that
\begin{equation} \nonumber
\int_{\R }F(|x|,w_1,\dots,w_m)dx >
\int_{\R}F(|x|,u_1,\dots,u_m)dx
\end{equation}
and thus
\begin{equation} \nonumber
J(w_1,w_2,\dots,w_m) < J(u_1,u_2,\dots,u_m) \leq T.
\end{equation}
Recalling that $( w_1,\dots, w_m)\in {\mathcal C}$ we have
proved that condition (H2) hold.\quad\qed
\vskip3pt

\subsection{Proof of Propositions~\ref{qualitative},~\ref{restrictions} and~\ref{restrictionsbis}}\label{section4}

First we state some known facts.

\begin{lemma} \label{stuart}
Assume that~\eqref{zerocvs}-\eqref{highdimension} hold. Then we
have
\begin{itemize}
\item[1)] Any minimizing sequence for~\eqref{0.11} is bounded in
$H^1(\R^N)$.
\item[2)] Any minimizing sequence satisfies {\rm (H1)}.
\item[3)] $m(d) < 0$ for any $d >0$.
\end{itemize}
\end{lemma}
\begin{proof}
The proof of these statements can be found in~\cite{stuarbbelow},
up to straightforward modifications at some places. We just
outline here the main steps. Also note that assertions 1) and 3)
are special cases of what we established in Step I of the proof of
Theorem~\ref{mainappl}. Assertion 1) is a direct consequence
of~\eqref{zerocvs} combined with standard H\"older and Sobolev
inequalities. Assertion 2) holds true because of the
limit~\eqref{compact} (see for instance~\cite[Lemma
5.2]{stuarbbelow} for such a result). Assertion 3) can be proved
using suitable test functions and taking advantage that,
under~\eqref{highdimension}, $F(x,s)$ does not decrease too fast
as $|x|$ goes to infinity (see~\cite[Theorem 5.4]{stuarbbelow}).
\end{proof}

The proof of Proposition~\ref{qualitative} relies on the following two lemmas.

\begin{lemma} \label{infimum-zero}
Assume that~\eqref{zerocvs}-\eqref{highdimension} hold.
Then, for any $d>0$, any $\varepsilon>0$
and all $R_0 >0$ there exists a function
$v \in C_0^{\infty}(\R^N)$ such that
$$
\|v\|_2^2 =d, \qquad {\rm supp}(v) \subset \R^N\setminus B(0, R_0),\qquad
 I(v) \leq \varepsilon.
$$
\end{lemma}
\begin{proof}
Take a positive function $u \in C_0^{\infty}(\R^N)$ such that $\|u\|_2^2 =d$. Then, considering the scaling
$t \mapsto t^{\frac{N}{2}}u(tx)=u_t(x)$, for all $t>0$, we get
$$
\int_{\R^N}|u_t|^2 dx = d, \qquad
\int_{\R^N}|\nabla u_t|^2 dx =  t^2 \int_{\R^N}|\nabla u|^2 dx.
$$
Since $\|u_t\|_{\infty} \to 0$ as $t \to 0^+$, given $\eps>0$, we can fix a value $t_0>0$ such that
$$
\frac{1}{2}\int_{\R^N}|\nabla u_{t_0}|^2 dx \leq \varepsilon
\qquad \text{and} \qquad \|u_{t_0}\|_{\infty} \leq \delta,
$$
where $\delta >0$ is the number which
appears in condition~\eqref{highdimension}. Translate now
$u_{t_0}$ into $\tilde{u}_{t_0}(\cdot) = u_{t_0}(\cdot+y)$ for a suitable
$y \in \R^N$ in such a way that
$$
{\rm supp}(\tilde{u}_{t_0}) \subset \R^N\setminus B(0, R_0).
$$
Then, since in view of~\eqref{highdimension}, $F(x,s) \geq 0$ for all $|x|$
sufficiently large and for $s \in [0, \delta]$, we obtain
$$
\int_{\R^N}F(x,u_{t_0})dx \geq 0.
$$
Thus
$$
I(\tilde{u}_{t_0}) \leq \frac{1}{2}\int_{\R^N}|\nabla
u_{t_0}|^2 dx \leq \varepsilon,
$$
and $v:= \tilde{u}_{t_0}$ has all
the desired properties.
\end{proof}

\begin{lemma}
    \label{invisible-mass}
Assume that~\eqref{zerocvs}-\eqref{highdimension} hold and let $u \in
C_0^{\infty}(\R^N)$ be such that $\|u\|_2^2 <c$. Then, for any
$\varepsilon >0$, there exists  a function $v \in C_0^{\infty}(\R^N)$ such that
$$
I(u+v) \leq I(u) +  \varepsilon, \qquad \|u + v\|_2^2 = c.
$$
\end{lemma}

\begin{proof}
Let $\varepsilon>0$ be fixed. By Lemma~\ref{infimum-zero}
we learn that there exists a function $v \in C_0^{\infty}(\R^N)$ with
$\|v\|_2^2 = c - \|u\|_2^2 >0$ and such that (since the supports of $u$ and $v$
can be assumed to be disjoint)
$$
\|u + v\|_2^2 = \|u\|_2^2 +
\|v\|_2^2 = c,
$$
and
$$
I(u + v) = I(u) + I(v) \leq I(u) + \varepsilon,
$$
which concludes the proof.
\end{proof}

We can now give the proof of Proposition~\ref{qualitative}.

\begin{proof}
We know by Lemma~\ref{stuart} that $m(c) <0$ for any $c >0$. Now,
assume by contradiction that there exist $0<c_1 < c_2$ such that
$m(c_1) < m(c_2)$ and set $m(c_2) - m(c_1) =\delta>0$. By
definition of $m(c_1)$  there exists a $u_{c_1} \in H^1(\R^N)$
such that $\|u_{c_1}\|_2^2 = c_1$ and $I(u_{c_1}) \leq m(c_1) +
\frac{\delta}{4}$. Arguing as in Lemma \ref{compact-support},
where we can directly use the continuity of the functional $I$, we
can assume that $u_{c_1} \in C_0^{\infty}(\R^N)$. Now, by
Lemma~\ref{invisible-mass}, since $\|u_{c_1}\|_2^2 < c_2$, we can
find a function $v \in C_0^{\infty}(\R^N)$ such that
$$
I(u_{c_1}+v) \leq I(u_{c_1}) + \frac{\delta}{4}
$$
and
$\|u_{c_1}+ v\|_2^2 = c_2$. Then we get that
$$
I(u_{c_1}+v) \leq m(c_1) + \frac{\delta}{2}< m(c_2).
$$
This contradiction proves Proposition \ref{qualitative}.
\end{proof}
\vskip2pt

We now give the proof of Proposition~\ref{restrictions}, which covers the case $N=1$.
\begin{proof}
Let $(u_n) \subset H^1(\R)$ be a positive minimizing sequence for
problem~\eqref{0.11}. This is possible by~\eqref{modineqs}. From
Lemma~\ref{stuart}, we can assume that $u_n \rightharpoonup u$
with $u\geq 0$ and $I(u) \leq m(c) <0$. To conclude, we need to
show that $\|u\|_2^2 = c$. Since $I(u) <0$, we have that $u \neq
0$. Thus assume by contradiction that $0 < \|u\|_2^2 < c$. We
distinguish two cases according to the fact that there exists, or
not, a point $x_0 \in \R$ such that $u(x_0) >0$ and $u$ is
non-increasing over $[x_0, +\infty[$. We also recall that elements
of $H^1(\R)$ are continuous functions which vanish as
$|x|\to\infty$.
\vskip3pt
\noindent {\bf Case I.} We assume that
there exists a $x_0 \in \R$ such that $u(x_0)>0$ and $u$ is
non-increasing over $[x_0, +\infty[$. In this situation we use the
same trick as in the proof of Theorem~\ref{mainappl}. Since
$u(x)\to 0$ as $|x|\to\infty$, without loss of generality, we may
assume that $u(x) \in [0, \delta]$, for all $x \in [x_0,+\infty[$.
Now we define a function $w :\R \to \R$ by
$$
w(x):=
\begin{cases}
    u(x)       & \text{if $x \in ]-\infty,x_0]$}, \\
     u(x_0)     & \text{if $x \in [x_0,x_0+ \mu]$},\\
    u(x-\mu)    & \text{if $x \in [x_0 + \mu ,+ \infty[$}.
\end{cases}
$$
Here $\mu >0$ is chosen in order to have $\|w\|_2^2 = c$. Clearly $\|w'\|_2^2
= \|u'\|_2^2 $ and since $w \geq u$ with $w \neq u$ taking into
account condition~\eqref{increasing2}, we have that
$$\int_{\R}F(x,w)dx > \int_{\R}F(x,u) dx.$$
Thus $I(w) < I(u)$ and, since $\|w\|_2^2 = c$, we have reached a
contradiction.
\vskip3pt
\noindent
{\bf Case II.} In this case
there is no point $x_0 \in \R$ such that $u(x_0) >0$ and $u$ is
non-increasing on $[x_0, +\infty[$. In this situation,
necessarily, the following occurs: there exists $x_1,x_2 \in
[x_0,+\infty[$ with $x_1 < x_2$ such that $u(x) < u(x_1) = u(x_2)$
for $x \in ]x_1,x_2[$. Now we define $w :\R \to \R$ by setting
$$
w(x):=
\begin{cases}
    u(x)     & \text{if $x \in [-\infty,x_1]$}, \\
     u(x_1)   & \text{if $x \in [x_1,x_2]$},  \\
    u(x)    & \text{if $x \in [x_2,+ \infty[$}.
\end{cases}
$$
Then $w \in H^1(\R)$ with
$$\int_{\R}|w'|^2 dx < \int_{\R}|u'|^2 dx$$
and also, by~\eqref{increasing2},
$$\int_{\R}F(x,w)dx > \int_{\R}F(x,u)dx.$$
Now observe that the points $x_1,x_2$ can be chosen such that
$$
\int_{[x_1,x_2]}|u(x_1)|^2  - |u(x)|^2 dx >0
$$
is  smaller than $c - \|u\|_2^2>0$. Then $I(w)<I(u)$ and $\|w\|_2^2=d<c$,
so that the conclusion follows by Proposition~\ref{qualitative}.
\end{proof}
\medskip

Before proving Proposition~\ref{restrictionsbis} we show, under
our additional regularity assumptions, that any minimizer
satisfies a Euler-Lagrange equation and we discuss the value of
the associated Lagrange parameter.

\begin{lemma}\label{lagrange}
Assume that $f(x,s) = F_s(x,s)$ exists and
that~\eqref{compact}-\eqref{modineqs} and~\eqref{zerocvsbis} hold.
Then $I \in C^1(H^1(\R^N), \R)$ and we have
\begin{enumerate}
\item[i)] Any minimizer $v \in H^1(\R^N)$ of $I$ on $\|v\|_2^2 =c$
satisfies
$$
- \Delta v - f(x,v) = \beta v, \quad \mbox{ with } \beta =
\frac{I'(v)v}{ \|v\|_2^2} \leq 0.
$$

\item[ii)] Let $(u_n) \subset H^1(\R^N)$ with $\|u_n\|_2^2 = c$ be
such that $u_n \rightharpoonup u$ with $I(u) \leq m(c) <0$ and $0
< \|u\|_2^2 < c.$ Then $u$ satisfies the equation
\begin{equation}\label{equationb}
- \Delta u - f(x,u) = 0.
\end{equation}
\end{enumerate}
\end{lemma}

\begin{proof}
Assuming that $f(x,s) = F_s(x,s)$ exists and under
\eqref{compact}-\eqref{modineqs} and~\eqref{zerocvsbis} it is
classical to show that $I$ is a $C^1$-functional
(see~\cite{stuarbbelow}). Thus, by standard considerations, any
minimizer of $I$ on the constraint $\|v\|_2^2 = c$ satisfies
\begin{equation}\label{equationc}
- \Delta v - f(x,v) = \beta v, \quad \mbox{where } \beta \mbox{ is
given by } \beta = \frac{I'(v)v}{ \|v\|_2^2}.
\end{equation}
Now assume by contradiction that $\beta
>0$. Then $I'(v)v = \beta \|v\|_2^2
>0$ and thus, since one has,
\begin{equation}
I((1-t)v) = m(c) - t(I'(v)v + o(1)) \quad \mbox{as } t \to 0,
\end{equation}
we can fix a small $t_0 >0$ such that $v_0 = (1-t_0) v$ satisfies
$I(v_0) < m(c)$. Since $\|v_0\|_2^2 < c$ we have a contradiction
with Proposition~\ref{qualitative} which says that $\lambda \to
m(\lambda)$ is non increasing. This proves i).
Now assume that the assumptions of ii) hold. By
Remark~\ref{crucial} the weak limit $u \in H^1(\R^N)$ minimizes
$I$ on the constraint $\|u\|_2^2 := d < c$ (and $m(d)= m(c)$).
Also, by Part i) we know that the associated Lagrange multiplier
$\beta \in \R$ satisfies $\beta \leq 0$. Let us proves that $\beta
< 0$ is impossible. If we assume, by contradiction, that $\beta
<0$ then $I'(u)u <0$ and since one has
\begin{equation}
I((1+t)u) = m(c) + t(I'(u)u + o(1)) \quad \mbox{as } t \to 0,
\end{equation}
we can fix a small $t_0 >0$ such that $u_0 = (1+t_0) v$ satisfies
both $I(u_0) < m(c)$ and $\|u_0\|_2^2 < c.$ Here again this
provides a contradiction with the fact that $\lambda \to
m(\lambda)$ is non increasing.
\end{proof}
 \medskip

We can now give the proof of Proposition~\ref{restrictionsbis}.

\begin{proof}
Let $(u_n) \subset H^1(\R^N)$ be a positive minimizing sequence
for~\eqref{0.11}. From Lemma~\ref{stuart} we can assume that $u_n
\rightharpoonup u$ with $u\geq 0$ and $I(u) \leq m(c) <0.$ To
conclude we need to show that $\|u\|_2^2 = c$. Since $I(u) < 0$ we
have $u \neq 0$. Thus assume by contradiction that $0 < \|u\|_2^2
< c$. In turn, from Part ii) of Lemma~\ref{lagrange}, we learn
that $u \in H^1(\R^N)$ satisfies equation~\eqref{equationb}. Also,
since $f(x,s)\geq 0$ for $s\in\R^+$, it follows from the strong
maximum principle that $u>0$. Therefore, taking into
account~\eqref{highdimensionbis}, we see that $u$ is a weak
solution of the variational inequality
$$
-\Delta u\geq b(x)u^{1+\alpha}\quad\text{in $\R^N$},
$$
where $b:\R^N\to\R^+$ is defined by
$$
b(x)
=
\begin{cases}
    \frac{f(x,u(x))}{u^{1+\alpha}(x)} & \text{if $|x|\leq r_0$}, \\
    A(1+|x|)^{-d} & \text{if $|x|\geq r_0$},
\end{cases}
$$
being $r_0,d$ and $\alpha$ the positive numbers appearing in~\eqref{highdimensionbis}.
Now, from the Liouville type theorem~\cite[Theorem 3.1, Chapter I]{MiPo}, we know
that $u\equiv 0$ under the restrictions
on the values of $\alpha$ given in condition~\eqref{highdimensionbis}
(notice that only the behaviour of $b(x)$ for large values of
$|x|$, and hence the behaviour of the weight $|x|^{-d}$, determines the validity of
the result from~\cite{MiPo} (see~\cite[formulas (3.4) and (3.5)]{MiPo}).
This immediately provides us a contradiction, since $u\not\equiv 0$.
\end{proof}

\begin{remark}
From our results of minimization we can derive bifurcation results for the
equation
\begin{equation}\label{bifurcation}
- \Delta u + \beta u = f(x,u), \quad u \in H^1(\R^N), \,\,\, \beta \in
\R.
\end{equation}
We recall that $\beta = 0$ is a bifurcation point for~\eqref{bifurcation} if
there exists a sequence $(\beta_n, u_n) \subset \R \times
H^1(\R^N)\backslash\{0\}$ of solutions of~\eqref{bifurcation} such that
$\beta_n \to 0$ and $\|u_n\|_{H^1(\R^N)} \to 0$ as $n\to\infty$.  The point here is that the
bifurcation phenomena occurs from the bottom of the essential spectrum.
\vskip2pt

Let $(c_n) \subset (0, + \infty)$ be such that $c_n \to 0$. Under the
assumptions that $f(x,s)$ exists and that~\eqref{compact}-\eqref{highdimension}
and~\eqref{zerocvsbis} hold we immediately derive from Remark~\ref{crucial} and
Part i) of Lemma~\ref{lagrange} the existence of a sequence $(\beta_n , u_n)
\subset [0, + \infty) \times H^1(\R^N)\backslash\{0\}$ such that $(\beta_n,
u_n)$ satisfies~\eqref{bifurcation} with $ 0 < \|u_n\|_2^2 \leq c_n$. From this
it is standard to show that $\beta_n \to 0$ and $\|u_n\|_{H^1(\R^N)}~\to 0$
as $n\to\infty$ (see~\cite{stuarbbelow}).
\vskip2pt

If instead of~\eqref{highdimension} we require
assumption~\eqref{highdimensionbis} we know, in addition, that $(\beta_n) \subset
(0, + \infty)$ and that $\|u_n\|_2^2 = c_n$. The fact that $\|u_n\|_2^2 = c_n$
follows directly from Proposition~\ref{restrictionsbis} and Part i) of
Lemma~\ref{lagrange}. To exclude the possibility that that $\beta_n =
0$ (thus showing that the bifurcation occurs by regular values) one can argue as in
the proof of Proposition~\ref{restrictionsbis}.
\vskip2pt

We also mention that, as long as we are interested only in the bifurcation
phenomena, we can remove the condition at infinity in~\eqref{zerocvsbis}.
Indeed observing that $\|u_n\|_{\infty} \to 0$ as $ n \to \infty$ we are free
to modify $f(x,s)$ outside the origin in $s \in \R$ (see~\cite{jeanjean} for
such arguments).
\end{remark}

\subsection{Proof of Theorem~\ref{Badiale-Rolando}}\label{section5}

We start with some preliminaries following closely~\cite{BaRo}. \medskip

We equip the Sobolev spaces $H$ and $H_s$ with the Hilbert norm
\begin{equation}\label{normBR}
\|u\| := \Big( \int_{\R^N} |\nabla u|^2 dx + \mu \int_{\R^N}
\frac{|u|^2}{|y|^2}dx + \int_{\R^N}|u|^2 dx \Big)^{\frac{1}{2}},
\quad \text{ for all $u \in H$.}
\end{equation}
Clearly $H_s \subset H \subset H^1(\R^N)$ and thus $H \subset
L^p(\R^N)$, for $2 \leq p \leq \frac{2N}{N-2}$. To simplify the
notation it is also useful to denote
$$
\|u\|_X := \Big( \int_{\R^N} |\nabla u|^2 dx + \mu \int_{\R^N}
\frac{|u|^2}{|y|^2}dx \Big)^{\frac{1}{2}}.
$$
Also observe that for any function $f \in C(\R, \R)$ satisfying $(f_1)$
and $(f_2)$ or $(f_3)$ we have
\begin{equation}\label{Bound-F}
|F(t)| \leq M ( |t|^p + |t|^q ), \quad \text{for all $t \in \R$}
\end{equation}
with $p,q \in ]2, 2 + \frac{4}{N}[$ and some positive constant $M$. Now
it is a standard fact, that under inequality~\eqref{Bound-F} the functional $J: H
\to \R$ defined by
$$
J(u) := \frac{1}{2}\|u\|_X^2 - \int_{\R^N} F(u) dx
$$
is well defined and continuous on $H$. Finally to study the minimization
problem~\eqref{minimization-problem}, for any $\rho >0$, we set
$$
\calm_{\rho}:= \Big\{ u \in H_s:\, \int_{\R^N}|u|^2 dx = \rho \Big\}
\quad \text{and} \quad m_{\rho}:= \inf_{u \in \calm_{\rho}}J(u).
$$
\noindent
We now give the proof of Theorem~\ref{Badiale-Rolando}. \medskip
\vskip2pt
\noindent
First from~\cite{BaRo} we borrow the next results, which hold
true under the assumptions of Theorem~\ref{Badiale-Rolando}.

\begin{lemma}\label{negative-infimumBD}
There exists a $\rho_0 >0$ such that $m_{\rho} <0$ for any $\rho >\rho_0$.
\end{lemma}

\begin{proof}
This follows directly from~\cite[Proposition 3.1 and Corollary 3.1]{BaRo}.
\end{proof}

The next result is exactly Lemma 4.2 of~\cite{BaRo}.

\begin{lemma}
    \label{Lemma4.2}
For every $\rho >0$, problem~\eqref{minimization-problem} admits
bounded minimizing sequences $(u_n)$ such that $u_n(y,z) =
u_n(|y|,|z|) \geq 0$ is non increasing in $|z|$.
Moreover, if any of such sequences satisfies
\begin{equation}
    \label{non-vanishingBR}
\inf_{n \in \N}\int_{B(x_n,R)} |u_n|^2 dx >0, \quad \text{ for some
$R >0$ and $(x_n) \subset \R^N$},
\end{equation}
then the sequence $(x_n)$ is bounded.
\end{lemma}

 Now we conclude the proof of
Theorem~\ref{Badiale-Rolando} with the following lemma.

\begin{lemma}\label{matBR}
Let $\rho >0$ be such that $m_{\rho}<0$ and $(u_n) \subset H_s$ be
a minimizing sequence as given by Lemma~\ref{Lemma4.2}. Then, up
to a subsequence, $u_n \rightharpoonup u$ with $J(u) \leq
m_{\rho}$ and $\|u\|_2^2 = \rho.$
\end{lemma}

\begin{proof}
Taking a minimizing sequence as given in Lemma~\ref{Lemma4.2}, we
can assume that $u_n \rightharpoonup u$ in $H_s$ as $n\to\infty$. Also, from the second
part of Lemma~\ref{Lemma4.2}, we see that, for any $\varepsilon>0$,
there exists a radius $R(\varepsilon) >0$ such that
\begin{equation}
    \label{vanishingBR}
\limsup_{n \to \infty} \sup_{x \in \R^N \backslash B(0,
R(\varepsilon))} \int_{B(x,1)}|u_n|^2 dx \leq \varepsilon.
\end{equation}
Following the proof of~\cite[Lemma I.1]{lions2}, we thus have
\begin{equation}
    \label{add}
\limsup_{n \to \infty} \int_{\R^N \backslash B(0,
R(\varepsilon))}|u_n|^p dx \leq C(\varepsilon),  \quad\, \text{for
any $2 <p < \frac{2N}{N-2}$},
\end{equation}
where $C(\varepsilon) \to 0$ provided that $\varepsilon \to 0$. Now, we fix an
arbitrary $\varepsilon >0$.  Because of the compact embedding $H \subset L_{{\rm loc}}^p (\R^N)$
for all $1\leq p < \frac{2N}{N-2}$, using~\eqref{Bound-F}, as $n\to\infty$ we obtain
\begin{equation}
    \label{strong-convergence}
\int_{B(0, R(\varepsilon))} F(u_n) dx \to \int_{B(0,R(\varepsilon))}F(u) dx.
\end{equation}
Gathering~\eqref{add} and~\eqref{strong-convergence},
since $\varepsilon >0$ is arbitrary, it follows that
$$
\int_{\R^N} F(u_n) dx \to \int_{\R^N}F(u) dx,
$$
as $n\to\infty$. Also, because $\|\cdot\|_X$ is a norm, $\|u\|_X^2 \leq \liminf_{n
\to \infty}\|u_n\|_X^2.$ Thus we do have
$$
J(u) \leq \liminf_{n \to \infty} J(u_n) = m_{\rho}.
$$
Namely (H1) hold. Now if $\|u\|_2^2 = \rho$ we are done.
Consequently we assume, by contradiction, that $\|u\|_2^2 < \rho$.
Since $J(u) \leq m_{\rho} <0$, $u=0$ is impossible. Thus $0 <
\|u\|_2^2 = \lambda$ and we consider the scaling $v(x) =
u(t^{-\frac{1}{N}}x)$ for $t>1$. Clearly for $ t = \frac{\rho}{\lambda} >1$ we have $\|v\|_2^2
= \rho.$ Now, since $t >1$ and $J(u)<0$,
\begin{align*}
        J(v)  & = \frac{1}{2}t^{1 - \frac{2}{N}}\|u\|_X^2 - t \int_{\R^N}F(u) \\
&       = t \Big[ \frac{1}{2}t^{- \frac{2}{N}} \|u\|_X^2 - \int_{\R^N}F(u)\Big]  < t J(u)<m_\rho.
    \end{align*}
Thus we reach a contradiction and the proof is complete.
\end{proof}

\bigskip

\end{document}